\documentclass[preprint]{elsarticle}

\usepackage{amsmath}
\usepackage{amsfonts}
\usepackage{amssymb}
\usepackage{eucal}

\newtheorem{cor}{Corollary}
\newtheorem{prop}{Proposition}
\newtheorem{thm}{Theorem}

\newdefinition{rmk}{Remark}
\newproof{pf}{Proof}
\newproof{pot}{Proof of Theorem \ref{thm2}}

\newcounter{ex}

\begin{document}

\author{Joanna We{\l}yczko}
\ead{joanna.welyczko@pwr.wroc.pl,
 Tel.:
(+48 71) 320-25-29, 
(+48 71) 320-22-18 
FAX: (+48 71) 328-07-51
}
\address{ Institute of Mathematics and Computer Science,
Wroc{\l}aw University of Technology,
Wybrze\.{z}e Wyspia\'nskiego 27, 50-370 Wroc{\l}aw, Poland.}
\title{On basic curvature identities for\\
almost (para)contact metric manifolds}

\begin{abstract}
In monograph of D. E. Blair 
{\it Riemannian geometry of contact and symplectic manifolds }
 and in the paper of S. Zamkovoy
{\it Canonical connections on paracontact manifolds, } 
 the curvature identities respectively for contact and paracontact metric manifold are proved. We obtain the curvature identity in the wider class of manifolds, which generalizes  results presented  in above mentioned publications. 
Moreover, we present some properties of almost (para)hermitian structure on a special semiproduct of $R_+$ and an almost (para)contact metric manifold. This semiproduct plays an auxiliary role in proving main theorem.
\end{abstract}

\maketitle
\section{Preliminaries}
Let $M$ be a $(2n+1)$-dimensional differentiable manifold endowed with  a $(1,1)$-tensor filed $\varphi$, a vector filed  $\xi$ and a 1-form  $\eta$ such that
\begin{equation}
\label{eps1}
  \varphi^2X =\varepsilon_1(X-\eta(X)\xi),\quad \eta(\xi)=1,
\end{equation}
where $\varepsilon_1=\pm1$. We note that (\ref{eps1}) implies $\varphi\xi=0$ and $\eta\circ\varphi=0$. When $\varepsilon_1=-1$, the triple $(\varphi,\xi,\eta)$ is an almost contact structure on $M$. When $\varepsilon_1=1$ and the tensor field $\varphi$ induces an almost paracomplex structure on the distribution $\mathcal{D}=\mathop{\rm Ker}\eta$, the triple is an almost paracontact structure on $M$ (\cite{CM,E,KW,SZ}). An almost paracpmplex structure on $D$ means that eigendistributions $\mathcal{D}^{\pm}$ corresponding to the eigenvalues $\pm1$ of $\varphi$ are both $n$-dimensional. Assume additionally that $M$ is endowed with a Riemannian or pseudo-Riemannian metric $g$ such that 
\begin{equation}
\label{eps2}
  g(\varphi X,\varphi Y)=-\varepsilon_1(g(X,Y)-\varepsilon_0\eta(X)\eta(Y)),
\end{equation}
where $\varepsilon_0=\pm1$. 
One claims that (\ref{eps2}) implies $\eta(X)=\varepsilon_0 g(X,\xi)$, and consequently $g(\xi,\xi)=\varepsilon_0$. For simplicity, let us call the quadruple $(\varphi,\xi,\eta,g)$ satisfying the conditions (\ref{eps1}) and (\ref{eps2}) to be an almost (para)contact metric structure on $M$, and the manifold endowed with such a structure to be an almost (para)contact metric manifold. The skew-symmetric $(0,2)$-tensor field $\varPhi$, defined by $\varPhi(X,Y)=g(X,\varphi Y)$, is called the fundamental form corresponding to the structure. 
 
Adopted by us the compatibility condition of the metric $g$ with the structure $(\varphi,\xi,\eta)$ is very general in nature and involves a number of classes of manifolds  found in the literature \cite{BCGH,DA,Dug,Mass,matsu,T}. 

On an almost (para)contact metric manifold, we define the tensor filed $N^{(1)}$ by 
\begin{equation*}
\label{nijen1}
  N^{(1)}(X,Y)=N(X,Y)-2\varepsilon_1\,d\eta(X,Y)\xi,
\end{equation*}
where $N$ is the Nijenhuis tensor of $\varphi$ given by 
$$
  N(X,Y)=\varphi^2[X,Y]+[\varphi X,\varphi Y]-\varphi[\varphi X,Y]-\varphi[X, \varphi Y].
$$ 
If $N^{(1)}$ vanishes identically, then the almost (para)contact metric manifold is said to be normal. In fact, the normality does not depend on the metric $g$. The normality condition says that the almost (para)complex structure $J$ defined on $M\times\mathbb R$ by 
$$
  J\partial_t=\varepsilon_1\xi,\quad 
  JX=\varphi X+\eta(X)\partial_t \ \ \mbox{for any} \ \ X\in\mathfrak{X}(M).
$$
is integrable, $t$ being the Cartesian coordinate on $\mathbb R$ and $\partial_t={\partial}/{\partial t}$.
An almost (para)contact metric manifold is said to be (para)contact metric one if $\varPhi=d\eta$ (\cite{blair2, blair,E,SZ});

\section{An almost (para)complex manifold}
In the sequel, we will use the theorem binding a covariant derivative of the tensor field $J$ in the direction of the Nijenhuis tensor of $J$ with the curvature of an almost (para)complex manifold.

Let $\widetilde{M}$ be a $2n$-dimensional almost (para)complex manifold i.e. differentiable manifold endowed with a $(1,1)$-tensor field $J$ such that
$$
  J^2=\varepsilon_1 I,\quad \varepsilon_1=\pm1.
$$
When $\varepsilon_1=-1$, $J$ is an almost complex structure. When $\varepsilon_1=1$ and the $\pm1$ eigendistributions of $J$ are $n$-dimensional, $J$ is an almost paracomplex structure. Let  $\widetilde{N}$ be the Nijenhuis tensor of $J$, 
\begin{equation*}
\label{nj}
  \widetilde{N}(\widetilde{X},\widetilde{Y})=
  J^2[\widetilde{X},\widetilde{Y}]+[J\widetilde{X},J\widetilde{Y}]-J[J\widetilde{X},\widetilde{Y}]-J[\widetilde{X},J\widetilde{Y}].
\end{equation*}

For a (symmetric) affine connection $\widetilde{\nabla}$ on $\widetilde{M}$, let 
$$
  \widetilde{R}(\widetilde{X},\widetilde{Y}) = 
  [\widetilde{\nabla}_{\widetilde{X}},\widetilde{\nabla}_{\widetilde{Y}}]-\widetilde{\nabla}_{[\widetilde{X},\widetilde{Y}]}
$$
be the curvature operator of $\widetilde{\nabla}$. Moreover, assume
$$
  [\widetilde{\nabla}_{\widetilde{X}},J]
  =\widetilde{\nabla}_{\widetilde{X}}\circ J - J\circ \widetilde{\nabla}_{\widetilde{X}}
  =\nabla_XJ.
$$
The theorem below is a generalization of A. Gray's theorem applying to the curvature of an almost hermitian structure (\cite{AG}).
\begin{thm}
\label{gray}
If an affine connection $\widetilde\nabla$ satisfies additionally the condition
\begin{equation}
\label{assum}
  [\widetilde{\nabla}_{J\widetilde{X}},J] = \delta J[\widetilde{\nabla}_{\widetilde{X}},J]
\end{equation}
for a constant $\delta=\pm1$, then 
\begin{equation}
\label{lemat1}
\begin{array}{l}
   [\widetilde{\nabla}_{\widetilde{N}(\widetilde{X},\widetilde{Y})},J] 
    = \null-\varepsilon_1[\widetilde{R}(\widetilde{X},\widetilde{Y}),J]
       -[\widetilde{R}(J\widetilde{X},J\widetilde{Y}),J] \\\hskip 3cm
      \null +\delta J[\widetilde{R}(J\widetilde{X},\widetilde{Y}),J]
       +\delta J[\widetilde{R}(\widetilde{X},J\widetilde{Y}),J]. 
\end{array}
\end{equation}
\end{thm}

\begin{pf}
We have in general
\begin{equation}
\label{bazowa}
\begin{array}{l}
 \null-\varepsilon_1[\widetilde{R}(\widetilde{X},\widetilde{Y}),J]
     -[\widetilde{R}(J\widetilde{X},J\widetilde{Y}),J] 
     +\delta J[\widetilde{R}(J\widetilde{X},\widetilde{Y}),J]\\
     \null+\delta J[\widetilde{R}(\widetilde{X},J\widetilde{Y}),J]
=\null-\varepsilon_1[[\widetilde{\nabla}_{\widetilde{X}},\widetilde{\nabla}_{\widetilde{Y}}],J]
    -[[\widetilde{\nabla}_{J{\widetilde{X}}},\widetilde{\nabla}_{J{\widetilde{Y}}}],J]\\\hspace{20mm}
    +\delta J[[\widetilde{\nabla}_{\widetilde{X}},\widetilde{\nabla}_{J{\widetilde{Y}}}],J]
 \null+\delta J[[\widetilde{\nabla}_J{\widetilde{X}},\widetilde{\nabla}_{\widetilde{Y}}],J]
    +[\widetilde{\nabla}_{\widetilde{N}({\widetilde{X}},{\widetilde{Y}})},J].
    \end{array}
\end{equation}
Rewrite the condition (\ref{assum}) in the following equivalent way
\begin{equation}
\label{zl}
  \widetilde{\nabla}_{J{\widetilde{X}}}\circ J=
  J\circ\widetilde{\nabla}_{J{\widetilde{X}}}
  +\delta J\circ\widetilde{\nabla}_{\widetilde{X}}\circ J
  -\delta\varepsilon_1\widetilde{\nabla}_{\widetilde{X}}.
\end{equation}
Now, using among others (\ref{zl}), we find the following relations after some long but easy calculations
\begin{eqnarray}
\label{aux1}
  -\varepsilon_1[[\widetilde{\nabla}_{\widetilde{X}},\widetilde{\nabla}_{\widetilde{Y}}],J]
  &=& \null-\varepsilon_1\widetilde{\nabla}_{\widetilde{X}}\circ\widetilde{\nabla}_{\widetilde{Y}}\circ J
      +\varepsilon_1\widetilde{\nabla}_{\widetilde{Y}}\circ\widetilde{\nabla}_{\widetilde{X}}\circ J \\
  & & \null+\varepsilon_1 J\circ\widetilde{\nabla}_{\widetilde{X}}\circ\widetilde{\nabla}_{\widetilde{Y}}
      -\varepsilon_1 J\circ\widetilde{\nabla}_{\widetilde{Y}}\circ\widetilde{\nabla}_{\widetilde{X}},\nonumber\\\nonumber\\
\label{aux2}
  -[[\widetilde{\nabla}_{J{\widetilde{X}}},\widetilde{\nabla}_{J{\widetilde{Y}}}],J]
  &=&-\delta J\circ\widetilde{\nabla}_{\widetilde{X}}\circ J\circ\widetilde{\nabla}_{J{\widetilde{Y}}}
     +\delta\varepsilon_1\widetilde{\nabla}_{\widetilde{X}}\circ\widetilde{\nabla}_{J{\widetilde{Y}}}\\
  & &\null-\delta J\circ\widetilde{\nabla}_{J{\widetilde{X}}}\circ\widetilde{\nabla}_{\widetilde{Y}}\circ J
     -J\widetilde{\nabla}_{\widetilde{X}}\circ J\circ\widetilde{\nabla}_{\widetilde{Y}}\circ J\nonumber\\
  & &\null+\varepsilon_1\widetilde{\nabla}_{\widetilde{X}}\circ\widetilde{\nabla}_{\widetilde{Y}}\circ J
     +\delta\varepsilon_1\widetilde{\nabla}_{J{\widetilde{X}}}\circ\widetilde{\nabla}_{\widetilde{Y}}\nonumber\\
  & &\null+\delta J\circ\widetilde{\nabla}_{\widetilde{Y}}\circ J\circ\widetilde{\nabla}_{J{\widetilde{X}}}
     -\delta\varepsilon_1\widetilde{\nabla}_{\widetilde{Y}}\circ\widetilde{\nabla}_{J{\widetilde{X}}}\nonumber\\
  & &\null+\delta J\circ\widetilde{\nabla}_{J{\widetilde{Y}}}\circ\widetilde{\nabla}_{{\widetilde{X}}}\circ J
     +J\circ\widetilde{\nabla}_{{\widetilde{Y}}}\circ J\circ\widetilde{\nabla}_{\widetilde{X}}\circ J\nonumber\\
  & &\null-\varepsilon_1\widetilde{\nabla}_{\widetilde{Y}}\circ\widetilde{\nabla}_{\widetilde{X}}\circ J
     -\delta\varepsilon_1\widetilde{\nabla}_{J{\widetilde{Y}}}\circ\widetilde{\nabla}_{\widetilde{X}},\nonumber\\\nonumber\\
\label{aux3}
  \delta J[[\widetilde{\nabla}_{\widetilde{X}},\widetilde{\nabla}_{J{\widetilde{Y}}}],J]                                
  &=&\delta J\circ\widetilde{\nabla}_{\widetilde{X}}\circ J\circ\widetilde{\nabla}_{J{\widetilde{Y}}}
     +J\circ\widetilde{\nabla}_{\widetilde{X}}\circ J\circ\widetilde{\nabla}_{\widetilde{Y}}\circ J\\
  & &\null-\varepsilon_1 J\circ\widetilde{\nabla}_{\widetilde{X}}\circ\widetilde{\nabla}_{\widetilde{Y}}
     -\delta J\circ\widetilde{\nabla}_{J{\widetilde{Y}}}\circ\widetilde{\nabla}_{{\widetilde{X}}}\circ J\nonumber\\
  & &\null-\delta\varepsilon_1\widetilde{\nabla}_{\widetilde{X}}\circ\widetilde{\nabla}_{J{\widetilde{Y}}}
     +\delta\varepsilon_1\widetilde{\nabla}_{J{\widetilde{Y}}}\circ\widetilde{\nabla}_{\widetilde{X}},\nonumber\\\nonumber\\
\label{aux4}
  \delta J[[\widetilde{\nabla}_J{\widetilde{X}},\widetilde{\nabla}_{\widetilde{Y}}],J]
  &=&\delta J\circ\widetilde{\nabla}_{J{\widetilde{X}}}\circ\widetilde{\nabla}_{{\widetilde{Y}}}\circ J
     -\delta J\circ\widetilde{\nabla}_{\widetilde{Y}}\circ J\circ\widetilde{\nabla}_{J{\widetilde{X}}}\\
  & &\null-J\circ\widetilde{\nabla}_{\widetilde{Y}}\circ J\circ\widetilde{\nabla}_{\widetilde{X}}\circ J
     +\varepsilon_1 J\circ\widetilde{\nabla}_{\widetilde{Y}}\circ\widetilde{\nabla}_{\widetilde{X}}\nonumber\\
  & &\null-\delta\varepsilon_1\widetilde{\nabla}_{J{\widetilde{X}}}\circ\widetilde{\nabla}_{\widetilde{Y}}
     +\delta\varepsilon_1\widetilde{\nabla}_{\widetilde{Y}}\circ\widetilde{\nabla}_{J{\widetilde{X}}}\nonumber
\end{eqnarray}
Applying the expressions (\ref{aux1}) - (\ref{aux4}) turns (\ref{bazowa}) into (\ref{lemat1}).
\end{pf}

\section{An almost (para)contact metric manifold}
We will propose the construction of an almost (para)hermitian structure on a special semiproduct of $R_+$ and an almost (para)contact metric manifold (for a special case where an almost contact metric manifold is used see (\cite{MM})).

Let $M$ be an almost (para)contact metric manifold and $(\varphi,\xi,\eta,g)$ its almost (para)contact metric structure. On the product manifold $\widetilde{M}=R_+\times M$, consider the cone metric (a kind of warped product metric) $\widetilde{g}$ defined by
\begin{equation}
\label{gtilde}
  \widetilde{g}=-\varepsilon_0\varepsilon_1 dt^2 + t^2 g
\end{equation}
Define a $(1,1)$-tensor field $J$ on $\widetilde{M}$ by assuming
\begin{equation}
\label{J}
  J\partial_t=-\frac{\varepsilon_0}{t}\xi,\quad 
  JX=\varphi X-\varepsilon_0\varepsilon_1t\eta(X)\partial_{t}\ \ \mbox{for any}\ \ X\in\mathfrak{X}(M),
\end{equation}
where $t$ is the Cartesian coordinate on $\mathbb R_+$ and $\partial_t={\partial}/{\partial t}$. Using (\ref{gtilde}) and (\ref{J}), one can easily check that the pair $(J,\widetilde{g})$ becomes an almost  (para-)\-Hermitian structure on $\widetilde{M}$ (precisely, almost Hermitian if $\varepsilon_1=-1$, and almost para-Hermitian if $\varepsilon_1=1$), that is, 
$$
  J^2=\varepsilon_1 I,\quad 
  \widetilde{g}(J\widetilde{X},J\widetilde{Y})=-\varepsilon_1\widetilde{g}(\widetilde{X},\widetilde{Y})
  \ \ \mbox{for any}\ \ \widetilde{X},\widetilde{Y}\in\mathfrak{X}(\widetilde{M}).
$$
Let $\varOmega$ be the fundamental form corresponding to the structure $(J,\widetilde g)$, that is, $\varOmega(\widetilde X,\widetilde Y)=\widetilde g(\widetilde X,J\widetilde Y)$. In view of (\ref{gtilde}) and (\ref{J}), we have
$$
  \varOmega(X,Y)=t^2\varPhi(X,Y),\quad \varOmega(X,\partial_t)=-t \eta(X)
$$
and hence
\begin{equation}
\label{omega}
  \varOmega=t^2\varPhi-2t\eta\wedge dt.
\end{equation}

It is a strighforward verification that the Levi-Civita connection $\widetilde{\nabla}$ of $\widetilde{g}$ is given by 
\begin{equation}
\label{pochtilde}
  \widetilde{\nabla}_{\partial_{t}}\partial_{t}=0, \quad
  \widetilde{\nabla}_X\partial_{t}=\widetilde{\nabla}_{\partial_{t}}X=\dfrac{1}{t} X,\quad
  \widetilde{\nabla}_XY=\nabla_XY+\varepsilon_0\varepsilon_1 t g(X,Y)\partial_{t}
\end{equation} 
for any $X,Y\in\mathfrak{X}(M)$, $\nabla$ being the Levi-Civita connection of the metric $g$. Using (\ref{J}) and (\ref{pochtilde}), we find the following formulas for the covariant derivative of $J$
\begin{eqnarray}
  (\widetilde{\nabla}_{\partial_t}J)\partial_t&=&0,\quad 
  (\widetilde{\nabla}_{\partial_t}J)X=0,\nonumber\\
\label{pochj}
  (\widetilde{\nabla}_XJ)\partial_t&=&-\dfrac{1}{t}(\varepsilon_0\nabla_X\xi+\varphi X),\\
  (\widetilde{\nabla}_XJ)Y&=&(\nabla_X\varphi)Y+\varepsilon_1 g(X,Y)\xi -\varepsilon_0\varepsilon_1\eta(Y)X\nonumber\\
                      &&\null-\varepsilon_0\varepsilon_1t((\nabla_X\eta)Y
                        - g(X,\varphi Y)))\partial_t.\nonumber
\end{eqnarray}

\begin{prop}
The structure $(J,\widetilde g)$ defined by (\ref{gtilde}) and (\ref{J}) is (para-)\-K\"ahler $(\widetilde\nabla J=0)$ if and only if the structure $(\varphi,\xi,\eta,g)$ satisfies the condition
\begin{equation}
\label{sasak}
  (\nabla_X\varphi)Y=-\varepsilon_1 g(X,Y)\xi +\varepsilon_0\varepsilon_1\eta(Y)X.
\end{equation}
\end{prop}

\begin{pf}
By (\ref{pochj}), we  see that ($\widetilde\nabla J=0$) impiles (\ref{sasak}). 

To have the converse implication, first we put $Y=\xi$ in (\ref{sasak}) and find $\nabla_X\xi=-\varepsilon_0\varphi X$ and $	(\nabla_X\eta)Y=g(X,\varphi Y).$
These equalities together with (\ref{pochj}) give $(\widetilde{\nabla}J)=0$.
\end{pf}

The manifold $M$ will be called (para-)Sasakian if it realizes the condition (\ref{sasak}).

\begin{prop}
The structure $(J,\widetilde g)$ defined by (\ref{gtilde}) and (\ref{J}) is almost (para-)\-K\"ahler ($d\varOmega=0$) if and only if the structure $(\varphi,\xi,\eta,g)$ is a (para)contact one ($  d\eta=\varPhi$).
\end{prop}

\begin{pf}
By (\ref{omega}), we have
\begin{eqnarray*}
	d\varOmega=r^2d\varPhi+2tdt\wedge(\varPhi-d\eta),
\end{eqnarray*}
which gives a thesis.
\end{pf}

\begin{prop}
\label{l}
The structure $(J,\widetilde g)$ defined by (\ref{gtilde}) and (\ref{J}) satisfies the condition
\begin{equation}
\label{lemat2J}
[\widetilde{\nabla}_{J\widetilde{X}},J]\widetilde{Y}=\delta J[\widetilde{\nabla}_{\widetilde{X}},J]\widetilde{Y}\quad  for \quad\delta=\pm1,
\end{equation}
 if and only if the almost (para)contact structure $(\varphi,\xi,\eta,g)$ on $M$ satisfies the equality
\begin{equation}
\label{lemat2}
\begin{array}{l}
(\nabla_{\varphi X}\varphi)Y-\delta\varphi(\nabla_X\varphi)Y-\delta\varepsilon_1(\nabla_X\eta)(Y)\xi\\[3pt] \hskip 2 cm
=(\delta-1)(\varepsilon_1 g(\varphi X,Y)\xi-\varepsilon_0\varepsilon_1\eta(Y)\varphi X).
\end{array}
\end{equation}
\end{prop}
\begin{pf}
The condition  (\ref{lemat2}) is fulfilled if and only if   $\widetilde{\nabla}_{J\widetilde{X}}J=\delta J\widetilde{\nabla}_{X}J$. In view of (\ref{J}) and (\ref{pochj}), it can be equivalently written as 
\begin{eqnarray}
\label{(1)}
(\widetilde{\nabla}_{\varphi X}J)Y&=&\delta J(\widetilde{\nabla}_XJ)Y,\\
\label{(2)}
(\widetilde{\nabla}_{\varphi X}J)\partial_t&=&\delta J(\widetilde{\nabla}_XJ)\partial_t,\\
\label{(3)}
(\widetilde{\nabla}_{\xi}J)Y&=&0,\\
\label{(4)}
(\widetilde{\nabla}_{\xi}J)\partial_t&=&0.
\end{eqnarray}

Next using once again (\ref{J}) and (\ref{pochj}), we calculate, that
(\ref{(1)}) is equivalent to the pair of equalities (\ref{lemat2}) and 
\begin{equation}
\label{p12}
\begin{array}{l}
(\nabla_{\varphi X}\eta)(Y)- g(\varphi X,\varphi Y)\\[3pt]\hskip 1cm
=\delta(\eta(\nabla_X\varphi)(Y)- g(\varphi X,\varphi Y)).
\end{array}
\end{equation}
The condition (\ref{(2)}) is equivalent to
\begin{equation}
\label{p2}
\varepsilon_0\nabla_{\varphi X}\xi+\varphi^2X=\delta(\varepsilon_0\varphi\nabla_X\xi+\varphi^2X).
\end{equation}
The condition (\ref{(3)}) is equivalent to
\begin{equation}
\label{p3}
(\nabla_\xi\varphi)Y=0,\quad {\rm and}\quad (\nabla_\xi\eta)(Y)=0.
\end{equation}
The condition (\ref{(4)}) is equivalent to
\begin{equation}
\label{p4}
\nabla_\xi\xi=0.
\end{equation}
Now suffice it to see that only (\ref{lemat2}) is a significance for the proof.
Exactly, putting $X=Y=\xi$ in (\ref{lemat2}) and using (\ref{eps1}), we get (\ref{p4}).
Putting $X=\xi$ in (\ref{lemat2}) and using (\ref{eps1}) and (\ref{p4}), we get (\ref{p3}).
Putting $Y=\xi$ in (\ref{lemat2}) and using (\ref{eps1}), we get (\ref{p2}).
Using (\ref{p2}), we calculate (\ref{p12}).
\end{pf}

$\widetilde{M}$ is a semiproduct manifold, so that it is easy to check that
\begin{prop}
\label{p}
For the curvature of the manifold $\widetilde{M}$, we have
\begin{equation}
\begin{array}{l}
\label{R}
\widetilde{R}(\partial_t,X)\partial_t=0,\quad \widetilde{R}(X,Y)\partial_t=0,\quad \widetilde{R}(\partial_t,X)Y=0,\\
\widetilde{R}(X,Y)Z=R(X,Y)Z+\varepsilon_0\varepsilon_1 g(Y,Z)X-\varepsilon_0\varepsilon_1 g(X,Z)Y,
\end{array}
\end{equation}
where $\widetilde{R}$ (respectively $R$) are the curvature tensors for $\widetilde{g}$ (respectively for $g$).
\end{prop}
As a consequence of the Proposition \ref{p} and (\ref{J}), we get 
\begin{prop}
The curvature of the manifold $\widetilde{M}$ satisfies
\begin{equation}
\begin{array}{l}
\label{RJ}
\widetilde{R}(\partial_t,X)(J\partial_t)=0,\quad \widetilde{R}(\partial_t,X)(JY)=0,\\
\widetilde{R}(X,Y)(J\partial_t)=-\cfrac{\varepsilon_0}{t}[R(X,Y)\xi+\varepsilon_1\eta(Y)X-\varepsilon_1\eta(X)Y],\\
\widetilde{R}(X,Y)(JZ)=R(X,Y)(\varphi Z)+\varepsilon_0\varepsilon_1 g(Y,\varphi Z)X-\varepsilon_0\varepsilon_1 g(X,\varphi Z)Y.
\end{array}
\end{equation}
\end{prop}
For  $\varepsilon_0=1$ and $\varepsilon_1=-1$, formulas (\ref{sasak}), (\ref{R}), (\ref{RJ}) are presented in (\cite{MM}).
\begin{prop}
The Nijenhuis tensor  $\widetilde{N}$ of the operator $J$ on the manifold $\widetilde{M}$ is given as
\begin{equation}
\begin{array}{l}
\label{njac}
\widetilde{N}(X,Y)=N^{(1)}(X,Y)-\varepsilon_0\varepsilon_1t N^{(2)}(X,Y)\partial_t,\\[3pt]
\widetilde{N}(\partial_t,Y)=-\cfrac{\varepsilon_0}{t}N^{(3)}Y+\varepsilon_1N^{(4)}Y\partial_t,
\end{array}
\end{equation}
where 
\begin{eqnarray*}
N^{(2)}(X,Y)&=&(\mathcal{L}_{\varphi X}\eta)(Y)-(\mathcal{L}_{\varphi Y}\eta)(X),\\
N^{(3)}Y&=&(\mathcal{L}_{\xi}\varphi)Y,\\
N^{(4)}Y&=&(\mathcal{L}_{\xi}\eta)(Y).
\end{eqnarray*}
\end{prop}
\begin{pf}
Let us recall that
\begin{eqnarray*}
\widetilde{N}(X,Y)&=&J^2[X,Y]+[JX,JY]-J[JX,Y]-J[X,JY]\\
\widetilde{N}(\partial_t,Y)&=&J^2[\partial_t,Y]+[J\partial_t,JY]-J[J\partial_t,Y]-J[\partial_t,JY].
\end{eqnarray*}
Using (\ref{J}), (\ref{pochtilde}) and (\ref{pochj}) in above expressions, we obtain a thesis.
\end{pf}
The following theorem contines the main curvature identity in the considered class of manifolds. Later, we will present its application in two subclasses of the manifold (a) a (para)contact metric and  (b) an almost normal (para)contact metric. 
\begin{thm}
\label{glowne}
The curvature of an almost (para)contact metric manifold satisfying the condition {\rm(\ref{lemat2})} fulfills the following identity 
\begin{eqnarray}
\label{a0}
&&(\delta\left[R(Z,\varphi X),\varphi\right]+\delta\left[R(\varphi Z,X),\varphi\right]
\null+\varepsilon_1\left[R(\varphi Z,\varphi X),\varphi\right]\varphi\\\nonumber
&&\null+\left[R(Z,X),\varphi\right]\varphi)Y
\null+\eta(Y)(R(\varphi Z,\varphi X)\xi+\varepsilon_1 R(Z,X)\xi)\\\nonumber
&&\hskip 1cm=\delta\varepsilon_1\{(\nabla_{\varphi N(Z,X)}\varphi)Y+\varepsilon_1g(\varphi N(Z,X),Y)\xi-\varepsilon_0\varepsilon_1\eta(Y)\varphi N(Z,X)\}\\\nonumber
&&\hskip 1.5cm\null-\varepsilon_0\varepsilon_1(\delta-1)\{2g(\varphi Z,Y)\varphi X-2g(\varphi X,Y)\varphi Z-g(\varphi Z,\varphi Y)X\\\nonumber
&&\hskip 1.5cm\null+g(\varphi X,\varphi Y)Z
\null+g(Z,Y)\varphi^2X-g(X,Y)\varphi^2Z\}.
\end{eqnarray}
\end{thm}
\begin{pf}
From the Proposition \ref{l}. we know that a condition (\ref{lemat2}) is equivalent to (\ref{lemat2J}). Hence in view of (\ref{lemat1}) for $Z,X\in\mathfrak{X}(M)$, we have
\begin{equation}
\label{a1}
\begin{array}{l}
\delta\varepsilon_1[\widetilde{\nabla}_{ J\widetilde{N}(Z,X)},J]
=\delta[\widetilde{R}(JZ,X),J]+\delta[\widetilde{R}(Z,JX),J]\\[3pt]
\hskip 3,5cm\null+\varepsilon_1[\widetilde{R}(JZ,JX),J]J+[\widetilde{R}(Z,X),J]J. 
\end{array}
\end{equation}
Applying (\ref{J}), (\ref{pochj}), (\ref{(3)}) and (\ref{njac}), we get 
\begin{equation*}
\begin{array}{l}
\label{right}
\delta\varepsilon_1[\widetilde{\nabla}_{J\widetilde{N}(Z,X)},J]Y=\delta\varepsilon_1(\widetilde{\nabla}_{J\widetilde{N}(Z,X)}J)Y=\delta\varepsilon_1 (\widetilde{\nabla}_{\varphi N(Z,X)+\varepsilon_1N^2(Z,X)\xi}J)Y\\[3pt]\hskip 3.5cm
=\delta\varepsilon_1\{(\nabla_{\varphi N(Z,X)}\varphi)Y+\varepsilon_1g(\varphi N(Z,X),Y)\xi\\[3pt]\hskip 4cm
\null-\varepsilon_0\varepsilon_1\eta(Y)\varphi N(Z,X)
\null-\varepsilon_0\varepsilon_1t((\nabla_{\varphi N(Z,X)}\eta)Y\\[3pt]\hskip 4cm
\null-g(\varphi N(Z,X),\varphi Y))\partial_t\}.
\end{array}
\end{equation*}
Moreover using (\ref{J}), (\ref{R}) and (\ref{RJ}), we obtain 
\begin{eqnarray*}
\delta[\widetilde{R}(JZ,X),J]Y&=&\delta\{\widetilde{R}(JZ,X)JY-J\widetilde{R}(JZ,X)Y\}\\
&=&\delta\{[R(\varphi Z,X),\varphi]Y+\varepsilon_0\varepsilon_1(g(\varphi Z,Y)\varphi X-g(\varphi X,Y)\varphi Z\\
&&\null-g(\varphi Z,\varphi Y)X-g(X,Y)\varphi^2Z)\\
&&\null+\varepsilon_0\varepsilon_1t(\eta(R(\varphi Z,X)Y)-\varepsilon_0\varepsilon_1g(\varphi Z,Y)\eta(X))\partial_t\},
\end{eqnarray*}
\begin{eqnarray*}
\delta[\widetilde{R}(Z,JX),J]Y&=&\delta\{\widetilde{R}(Z,JX)JY-J\widetilde{R}(Z,JX)Y\}\\
&=&\delta\{[R(Z,\varphi X),\varphi]Y+\varepsilon_0\varepsilon_1(g(\varphi Z,Y)\varphi X-g(\varphi X,Y)\varphi Z\\
&&\null+g(\varphi X,\varphi Y)Z+g(Z,Y)\varphi^2X)\\
&&\null+\varepsilon_0\varepsilon_1t(\eta(R( Z,\varphi X)Y)+\varepsilon_0\varepsilon_1g(\varphi X,Y)\eta(Z))\partial_t\},
\end{eqnarray*}
\begin{eqnarray*}
\varepsilon_1[\widetilde{R}(JZ,JX),J]JY&=&\varepsilon_1\widetilde{R}[JZ,JX]J^2Y-\varepsilon_1J\widetilde{R}(JZ,JX)JY\\
&=&\varepsilon_1[R(\varphi Z,\varphi X),\varphi]\varphi Y-\varepsilon_0\varepsilon_1(g(\varphi Z,Y)\varphi X\\
&&\null-g(\varphi X,Y)\varphi Z
\null-g(\varphi Z,\varphi Y)X\\
&&\null+g(\varphi X,\varphi Y)Z)
+\varepsilon_0g(Z,Y)\eta(X)\xi\\
&&\null-\varepsilon_0g(X,Y)\eta(Z)\xi\\
&&\null+\eta(Y)R(\varphi Z,\varphi X)\xi+\varepsilon_0t\eta(R(\varphi Z,\varphi X)\varphi Y)\partial_t,
\end{eqnarray*}
\begin{eqnarray*}
[\widetilde{R}(Z,X),J]JY&=&\widetilde{R}(Z,X)J^2Y-J\widetilde{R}(Z,X)JY\\
&=&[R(Z,X),\varphi]\varphi Y-\varepsilon_0\varepsilon_1(g(\varphi Z,Y)\varphi X-g(\varphi X,Y)\varphi Z\\
&&\null+g(Z,Y)\varphi^2X-g(X,Y)\varphi^2Z)-\varepsilon_0g(Z,Y)\eta(X)\xi\\
&&\null+\varepsilon_0g(X,Y)\eta(Z)\xi+\varepsilon_1\eta(Y)R(Z,X)\xi\\
&&\null+\varepsilon_0\varepsilon_1t(\eta(R(Z,X)\varphi Y)+\varepsilon_0\varepsilon_1g(\varphi Z, Y)\eta(X)\\
&&\null-\varepsilon_0\varepsilon_1g(\varphi X, Y)\eta(Z))\partial_t.
\end{eqnarray*}
Now putting  above commutators to (\ref{a1}) and comparing  parts tangent to $M$ in   (\ref{a1}), we get thesis. 
\end{pf}
\section{A (para)contact metric manifold}
We recall some properties of (para)contact metric manifolds.

Let $M$ be a (para)contact metric manifold. Let us define $h=\frac12\mathcal{L}_\xi\varphi$, where $\mathcal{L}$ is the Lie derivative, then 
\begin{equation}
\label{h1}
\begin{array}{lll}
  & g(hX,Y)=g(hY,X)\quad\mbox{($h$ is a symmetric operator)}, & \\
  & \varphi h+h\varphi=0,\quad \mathop{\rm Tr} h=0,\quad h\xi=0,\quad \eta\circ h=0. 
  \end{array}
\end{equation}

Moreover, the following conditions are fulfilled on $M$ 
\begin{eqnarray}
\label{h}
  \nabla_X\xi&=&\null-\varepsilon_0\varphi X+\varepsilon_1\varphi hX, \\
\label{lemat}
  \varepsilon_1(\nabla_{\varphi X}\varphi)\varphi Y-(\nabla_X\varphi)Y &=&
2\varepsilon_1g(X,Y)\xi\\\nonumber
&&\null-\varepsilon_1\eta(Y)(\varepsilon_0X-\varepsilon_1hX+\varepsilon_0\eta(X)\xi). 
\end{eqnarray}
Additionally $\nabla_\xi\xi=0$, $\nabla_\xi\varphi=0$, $\eta\circ h=0$ and we know that 
a (para)contact metric manifold is a (para)-Sasakian one if and only if  
\begin{equation}
\label{sas}
(\nabla_X\varphi)Y=-\varepsilon_1g(X,Y)\xi+\varepsilon_0\varepsilon_1\eta(Y)X.
\end{equation}

For a contact metric manifold ($\varepsilon_0=1$ i $\varepsilon_1=-1$), formulas (\ref{h1}), (\ref{h}), (\ref{lemat}) and (\ref{sas}) are proved in \cite{blair}. For a paracontact metric manifold ($\varepsilon_0=1$ i $\varepsilon_1=1$ ) they are shown in \cite{SZ}. In general, proofs are analogical to this contained in  \cite{blair} and \cite{SZ}.

\begin{prop}
For a (para)contact metric manifold, we have  
\begin{equation}
\label{s3}
\begin{array}{l}
(\nabla_{\varphi X}\varphi)Y+\varphi(\nabla_X\varphi)Y\\\hskip 2cm
=-\varepsilon_0\varepsilon_1(g(\varphi(\varepsilon_0X+\varepsilon_1hX),Y)\xi-2\eta(Y)\varphi X).
\end{array}
\end{equation}
\end{prop}
\begin{pf}
Using (\ref{eps1}), we see that (\ref{lemat}) is equivalent to the following condition
\begin{equation*}
\begin{array}{l}
(\nabla_X\varphi)Y+\varepsilon_1\varphi(\nabla_{\varphi X}\varphi)Y-\eta(\nabla_{\varphi X}Y)\xi+\varphi X\eta(Y)\xi+\eta(Y)\nabla_{\varphi X}\xi\\[3pt]\hskip 1cm
=\null-2\varepsilon_1g(X,Y)\xi+\varepsilon_1\eta(Y)(\varepsilon_0X-\varepsilon_1hX+\varepsilon_0\eta(X)\xi).
\end{array}
\end{equation*}
Putting $\varphi$ on the above equality and using (\ref{eps1}), (\ref{h}) and (\ref{h1}), we get (\ref{s3}).
\end{pf}

\begin{prop}
\label{Nijen}
The  Nijenhuis tensor $N$ of a (para)contact metric manifold satisfies  
\begin{eqnarray}
\label{s4}
N(X,Y)&=&-2\varphi(\nabla_X\varphi)Y+2\varphi(\nabla_Y\varphi)X+2\varepsilon_0\varepsilon_1\eta(Y)\varphi X\\\nonumber
&&-2\varepsilon_0\varepsilon_1\eta(X)\varphi Y
+2\varepsilon_1g(X,\varphi Y)\xi.
\end{eqnarray}
\end{prop}
\begin{pf}
With the help of the Levi-Civita connection on $M$ we get
$$
N(X,Y)=-\varphi(\nabla_X\varphi)Y+\varphi(\nabla_Y\varphi)X+(\nabla_{\varphi X}\varphi)Y-(\nabla_{\varphi Y}\varphi)X.
$$
Next, applying (\ref{s3}) and (\ref{h1}), we obtain (\ref{s4}).
\end{pf}
\begin{thm}
\label{glownec}
The curvature operator of a (para)contact metric manifold satisfies the following identity
\begin{eqnarray}
\label{w0}
&&(\left[R(Z,\varphi X),\varphi\right]+\left[R(\varphi Z,X),\varphi\right]
\null-\varepsilon_1\left[R(\varphi Z,\varphi X),\varphi\right]\varphi\\\nonumber
&&-\left[R(Z,X),\varphi\right]\varphi)Y
\null-\eta(Y)(\varepsilon_1 R(Z,X)\xi+R(\varphi Z,\varphi X)\xi\\\nonumber
&&\hskip1cm=\null-2(\nabla_{(\nabla_Z\varphi)X-(\nabla_X\varphi)Z}\varphi)Y+2\varepsilon_0\varepsilon_1\eta(X)(\nabla_Z\varphi)Y\\\nonumber
&&\hskip1.5cm\null-2\varepsilon_0\varepsilon_1\eta(Z)(\nabla_X\varphi)Y-2\varepsilon_1  g(Y,(\nabla_Z\varphi)X-(\nabla_X\varphi)Z)\xi\\\nonumber
&&\hskip1.5cm\null+2\varepsilon_0\varepsilon_1\eta(Y)((\nabla_Z\varphi)X-(\nabla_X\varphi)Z))\\\nonumber
&&\hskip1.5cm\null-4\varepsilon_0g(\varphi X,\varphi Y)\varphi^2Z+4\varepsilon_0g(\varphi Z,\varphi Y)\varphi^2X\\\nonumber
&&\hskip1.5cm\null+4\varepsilon_0\varepsilon_1 g(Y,\varphi X)\varphi Z-4\varepsilon_0\varepsilon_1 g(Y,\varphi Z)\varphi X.
\end{eqnarray}\end{thm}
\begin{pf}
On the (para)contact metric manifold, we have (\ref{s3}). It means, that the condition (\ref{lemat2}) is fulfilled with   $\delta =-1$. 
Moreover from (\ref{s4}), we get
$$
\varphi N(Z,X)=-2\varepsilon_1((\nabla_Z\varphi)X-(\nabla_X\varphi)Z)+2\varepsilon_0\eta(X)Z-2\varepsilon_0\eta(Z)X.
$$
In view of the above equality, we obtain 
\begin{eqnarray}
\label{rc}
  &&-\varepsilon_1\{(\nabla_{\varphi N(Z,X)}\varphi)Y
    +\varepsilon_1g(\varphi N(Z,X),Y)\xi
    -\varepsilon_0\varepsilon_1\eta(Y)\varphi N(Z,X))\}\\\nonumber
&&\hskip1.2cm=2(\nabla_{(\nabla_Z\varphi)X-(\nabla_X\varphi)Z}\varphi)Y-2\varepsilon_0\varepsilon_1\eta(X)(\nabla_Z\varphi)Y\\\nonumber 
&&\hskip1.5cm\null+2\varepsilon_0\varepsilon_1\eta(Z)(\nabla_X\varphi)Y+2\varepsilon_1 g((\nabla_Z\varphi)X-(\nabla_X\varphi)Z,Y)\xi\\\nonumber
&&\hskip1.5cm\null-2\varepsilon_0\eta(X)g(Z,Y)\xi+2\varepsilon_0\eta(Z)g(X,Y)\xi\\\nonumber
&&\hskip1.5cm\null-2\eta(Y)(\varepsilon_0\varepsilon_1((\nabla_Z\varphi)X-(\nabla_X\varphi)Z)-\eta(X)Z+\eta(Z)X).
\end{eqnarray}
Finally, the condition (\ref{w0})  follows from (\ref{a0}) after using (\ref{rc}), (\ref{eps1}), (\ref{eps2}) and putting $\delta=-1$.
\end{pf}
\begin{thm}
\label{ogolne}
The Riemannian curvature of the (para)contact metric manifold  satisfies the following relation
\begin{eqnarray}
\label{rcw2}
&&R(Z,\varphi X,\varphi Y,W)+ R(Z,\varphi X,Y,\varphi W)+R(\varphi Z,X,\varphi Y,W)\\\nonumber
&&\null+R(\varphi Z,X,Y,\varphi W)-R(\varphi Z,\varphi X,Y,W)-\varepsilon_1 R(\varphi Z,\varphi  X,\varphi Y,\varphi W)\\\nonumber
&&\null-\varepsilon_1 R(Z,X,Y,W)-R(Z,X,\varphi Y,\varphi W)\\\nonumber
&&\hskip 1cm=\null-2\sum_{i=1}^{2n+1}\epsilon_i(\nabla_{E_i}\Phi)(Z,X)(\nabla_{E_i}\Phi)(W,Y)\\\nonumber
&&\hskip1.5cm\null+2\varepsilon_0\varepsilon_1(\nabla_Z\Phi)(W,Y)\eta(X)-2\varepsilon_0\varepsilon_1(\nabla_X\Phi)(W,Y)\eta(Z)\\\nonumber
&&\hskip 1.5cm\null-2\varepsilon_0\varepsilon_1(\nabla_Y\Phi)(Z,X)\eta(W)
\null+2\varepsilon_0\varepsilon_1(\nabla_W\Phi)(Z,X)\eta(Y)\\\nonumber
&&\hskip 1.5cm\null+4\varepsilon_0g(\varphi X,\varphi Y)g(\varphi Z,\varphi W)
\null-4\varepsilon_0g(\varphi Z,\varphi Y)g(\varphi X,\varphi W)\\\nonumber
&&\hskip 1.5cm\null+4\varepsilon_0\varepsilon_1g(Y,\varphi X)g(\varphi Z,W)
\null-4\varepsilon_0\varepsilon_1g(Y,\varphi Z)g(\varphi X,W),
\end{eqnarray}
where $(E_i)$ is an orthonormal frame and $\varepsilon_i=g(E_i, E_i)$.
\end{thm}

\begin{pf}
Let us rewrite the condition (\ref{w0}) without commutators
\begin{eqnarray}
\label{rcw}
&&R(Z,\varphi X)\varphi Y-\varphi R(Z,\varphi X)Y+R(\varphi Z,X)\varphi Y-\varphi R(\varphi Z,X)Y\\\nonumber
&&\null-R(\varphi Z,\varphi X)Y
+\varepsilon_1\varphi R(\varphi Z,\varphi  X)\varphi Y-\varepsilon_1 R(Z,X)Y+\varphi R(Z,X)\varphi Y\\\nonumber
&&\hskip1cm=\null-2(\nabla_{(\nabla_Z\varphi)X-(\nabla_X\varphi)Z}\varphi)Y+2\varepsilon_0\varepsilon_1\eta(X)(\nabla_Z\varphi)Y\\\nonumber
&&\hskip 1.5cm\null-2\varepsilon_0\varepsilon_1\eta(Z)(\nabla_X\varphi)Y-2\varepsilon_1  g(Y,(\nabla_Z\varphi)X-(\nabla_X\varphi)Z)\xi\\\nonumber
&&\hskip 1.5cm\null+2\varepsilon_0\varepsilon_1\eta(Y)((\nabla_Z\varphi)X-(\nabla_X\varphi)Z))\\\nonumber
&&\hskip 1.5cm\null-4\varepsilon_0g(\varphi X,\varphi Y)\varphi^2Z+4\varepsilon_0g(\varphi Z,\varphi Y)\varphi^2X\\\nonumber
&&\hskip 1.5cm\null+4\varepsilon_0\varepsilon_1 g(Y,\varphi X)\varphi Z-4\varepsilon_0\varepsilon_1 g(Y,\varphi Z)\varphi X.
\end{eqnarray}
Let  $(E_i)$ be an orthonormal frame. Since $g((\nabla_X\varphi)Y,Z)=(\nabla_X\varPhi)(Z,Y)$ and 
$$
 (\nabla_X\varPhi)(Z,Y)+(\nabla_Y\varPhi)(X,Z)+(\nabla_Z\varPhi)(Y,X)=0,
$$
we have 
\begin{eqnarray*}
  g(Y,(\nabla_Z\varphi)X-(\nabla_X\varphi)Z)&=&(\nabla_Y\varPhi)(Z,X), \\
  g((\nabla_{(\nabla_Z\varphi)X-(\nabla_X\varphi)Z}\varphi)Y,W)
  &=&\sum_{i=1}^{2n+1}\epsilon_i(\nabla_{E_i}\Phi)(Z,X)(\nabla_{E_i}\Phi)(W,Y).
\end{eqnarray*}
Projecting (\ref{rcw}) on the vector field $W$ and using the above conditions, we get the thesis. 
\end{pf}
The Theorem \ref{ogolne}. is the most general relation concerning the Rimannian curvature of the (para)contact metric manifolds.  

\newpage

Up to now, It were known the identities for Ricci, *-Ricci, scalar, *-scalar curvature and some special properties of the Riemannian curvature tensor and only in two the following cases:\\[+3pt]
(a) for a contact metric manifold (i.e. $\varepsilon_0=1$ i $\varepsilon_1=-1$): Proposition 7.1, Corollary 7.1, Lemma 7.4 and Proposition 7.7 in the monograph \cite{blair}; \\
(b) for a paracontact metric manifold (i.e. $\varepsilon_0=1$ i $\varepsilon_1=1$): Proposition 3.1, Corollary 3.2, Lemma 3.4, Lemma 3.7 and Corollary 3.9 in the paper \cite{SZ}.\\[+3pt]
All these identities follow from Theorem \ref{ogolne}. Below, in Corollaries \ref{wn1} - \ref{wn4}, we present the form of mentioned identities in general case i.e. for any value of  $\varepsilon_0$ and $\varepsilon_1$.

\begin{cor}
\label{wn1}
The curvature  of a (para)contact metric manifold satisfies
\begin{equation*}
\label{rksi}
(R(\xi, X)\xi+\varepsilon_1\varphi R(\xi,\varphi X)\xi)=2\varphi^2X-2\varepsilon_1 h^2X,
\end{equation*}
\begin{equation*}
\label{rksi2}
\begin{array}{l}
-\varepsilon_1 R(\xi, X,Y,Z)- R(\xi, X,\varphi Y,\varphi Z)+R(\xi,\varphi X,\varphi Y,Z)\\[3pt]\hskip 1cm
\null+ R(\xi, \varphi X, Y,\varphi Z)
=2(\nabla_{hX}\Phi)(Y,Z)\\[3pt]\hskip 1cm
\null-2\varepsilon_0g(\varepsilon_0X-\varepsilon_1hX,Z)\eta(Y)+2\varepsilon_0g(\varepsilon_0X-\varepsilon_1hX,Y)\eta(Z).
\end{array}
\end{equation*}

\end{cor}
\begin{cor}
\label{wn2}
For any (para)contact metric manifold, we have
\begin{equation*}
\begin{array}{l}
Ric(\varphi X,\varphi Y)-\varepsilon_1 Ric(X,Y)+Ric^*(X,Y)+Ric^*(Y,X)\\[3pt]\hskip 1cm
=-\sum_{i=1}^{2n+1}\varepsilon_ig((\nabla_{E_i}\varphi)X,(\nabla_{E_i}\varphi)Y)+(4n-1)\varepsilon_0 g(X,Y)\\[3pt]\hskip 1,5cm
+\eta(X)\eta(Y)-2\varepsilon_1 g(X,hY)-\varepsilon_0g(hX,hY),
\end{array}
\end{equation*}
where $(E_i)$ is an orthonormal frame and $\varepsilon_i=g(E_i, E_i)$.
\end{cor}
\begin{cor}
\label{wn3}
On a (para)contact metric manifold the Ricci curvature in the direction of $\xi$ is given by
\begin{equation*}
Ric(\xi,\xi)=-\varepsilon_1(2n-|h|^2).
\end{equation*}
\end{cor}
\begin{cor}
\label{wn4}
On a (para)contact metric manifold the scalar curvatures  $r$ and $r^*$ fulfill the equality
$$
  {r}^*+{\varepsilon}_1r +4n^2=\mathop{\rm{Tr}}h^2+\frac{1}{2}(\mathop{\rm{Tr}}(\nabla\varphi)^2-4n).
$$
where $r$ is the scalar curvature and  $r^*=\mathop{\rm{Tr}}_{g}\{(X,Y)\rightarrow\mathop{Ric}(X,Y)\}$ is the  $*-$scalar curvature.
\end{cor}
From Theorem \ref{ogolne}. follows also  

\begin{cor}
On the (para)contact metric manifold of constant sectional curvature, we have  
\begin{equation}
\label{trp}
  \sum_{i=1}^{2n+1}\epsilon_iP(E_i,Z,X)P(E_i,W,Y)=0,
\end{equation}
where $P$ is the tensor field defined as follows 
$$
  P(X,Y,Z)=(\nabla_X\varPhi)(Y,Z)
           +\varepsilon_0\varepsilon_1g(X,Z)\eta(Y)
           -\varepsilon_0\varepsilon_1g(X,Y)\eta(Z),
$$
$(E_i)$ is the orthonormal frame and $\epsilon_i=g(E_i,E_i)$.
\end{cor}

\begin{pf}
Similarly, as in Theorem 3.12 in the paper \cite{SZ} it is proved that if the (para)contact metric manifold is of constant sectional curvature  $k$, then $k=-\varepsilon_0\varepsilon_1$. Next, putting in (\ref{rcw2}) the condition  
$$
R(X,Y,Z,W)=-\varepsilon_0\varepsilon_1(g(Y,Z)g(X,W)-g(X,Z)g(Y,W))
$$
and using (\ref{eps2}), we obtain 
\begin{eqnarray*}
&&\sum_{i=1}^{2n+1}\epsilon_i(\nabla_{E_i}\varPhi)(Z,X)(\nabla_{E_i}\varPhi)(W,Y)\\\nonumber
&&\null+\varepsilon_0\varepsilon_1\eta(W)(\nabla_Y\varPhi)(Z,X)-\varepsilon_0\varepsilon_1\eta(Y)(\nabla_W\varPhi)(Z,X)\\\nonumber
&&\null-\varepsilon_0\varepsilon_1\eta(X)(\nabla_Z\varPhi)(W,Y)+\varepsilon_0\varepsilon_1\eta(Z)(\nabla_X\varPhi)(W,Y)\\\nonumber
&&\null-g(Z,Y)\eta(W)\eta(X)+g(Z,W)\eta(Y)\eta(X)\\\nonumber
&&\null+g(X,Y)\eta(Z)\eta(W)-g(X,W)\eta(Z)\eta(Y)=0,
\end{eqnarray*}
that is equivalent to (\ref{trp}).
\end{pf}

\section{A normal almost(para)contact manifold}

In the proof of the theorem about the curvature operator, we will use the following proposition.
\begin{prop}
The following conditions are mutually equivalent
\begin{enumerate}
	\item[(a)] $M$ is a normal almost (para)contact metric manifold;
	\item[(b)] $\widetilde{M}$ is a (para)complex manifold;
	\item[(c)] \label{njn}
             $\widetilde{N}(X,Y)=0.$
             
\end{enumerate}
\end{prop}
\begin{pf}
If (a) is fulfilled then the tensor $N^{(1)}=0$. Hence, in view of \cite[Theorem 6.1]{blair} and \cite[Proposition 2.3]{SZ}, we claim that $N^{(2)}=0$, $N^{(3)}=0$ and $N^{(4)}=0$. Therefore from (\ref{njac}), we get $\widetilde{N}=0$. Hence $\widetilde{M}$ is a (para)complex  manifold and (b) is fulfilled. It is obvious that (b) implies (c).
Now, let us assume (c). Using (\ref{njac}), we get $N^{(1)}=0$, so that $M$ is a normal almost (para)contact metric manifold and hence (a) is satisfied.
\end{pf}

We also need the following necessary and sufficient condition for the normality.
For $\varepsilon_0=1$ and $\varepsilon_1=-1$, it is proved in the paper  S. Tanno \cite{ST}, and for $\varepsilon_0=1$ and $\varepsilon_1=1$ in the paper J. We{\l}yczko \cite{JW}. In general case i.e. for any value of  $\varepsilon_0$ and $\varepsilon_1$, the proof is analogous.

\begin{prop}
\label{norm}
The almost (para)contact metric manifolds is normal if and only if 
\begin{equation}
\label{normal}
\varphi(\nabla_X\varphi)Y-(\nabla_{\varphi X}\varphi)Y+\varepsilon_1(\nabla_X\eta)(Y)\xi=0,
\end{equation}
where $\nabla$ is the Levi-Civita'y connection. 
\end{prop}
\begin{thm}
The curvature operator of a normal almost (para)contact metric manifold satisfies the following identity
\begin{equation*}
\label{n3}
\begin{array}{l}
([R(Z,\varphi X),\varphi]+[R(\varphi Z, X),\varphi]+\varepsilon_1 [R(\varphi Z,\varphi X),\varphi]\varphi +[R(Z,X),\varphi]\varphi) Y\\[3pt]\hskip 1cm
=-\eta(Y)(\varepsilon_1 R(Z,X)\xi+R(\varphi Z,\varphi X)\xi).
\end{array}
\end{equation*}
\end{thm}
\begin{pf}
The normality of the almost (para)contact metric  is equivalent to con\-di\-tion (\ref{normal}). 
It means that (\ref{lemat2}) is fulfield with $\delta=1$.
Now, putting $\delta=1$ and $\varphi N(Z,X)=\varphi N^{(1)}(Z,X)=0$ in (\ref{a0}), we get the thesis. 
\end{pf}
\begin{thm}
\label{opkrzywnorm}
The curvature operator of the normal almost (para)contact metric manifold satisfies the identity 
\begin{eqnarray*}
  && \hskip -0.5 cm\varepsilon_1 R(Z,X)Y-\varphi R(Z,X)\varphi Y
      +R(Z,\varphi X)\varphi Y-\varphi R(Z,\varphi X)Y\\\nonumber
  && \hskip -0.5 cm\null+R(\varphi Z,X)\varphi Y-\varphi R(\varphi Z,X)Y
     +R(\varphi Z,\varphi X)Y-\varepsilon_1\varphi R(\varphi Z,\varphi X)\varphi Y=0.
\end{eqnarray*}
\end{thm}

Directly from the above theorem, for  $Z=\xi$, we get:

\begin{cor}
On the normal almost (para)contact metric manifold, we have 
\begin{eqnarray}
&&\varepsilon_1R(\xi,X)Y-\varphi R(\xi,X)\varphi Y
  +R(\xi,\varphi X)\varphi Y-\varphi R(\xi,\varphi X)Y=0,\nonumber\\
\label{nr9}
&&\varepsilon_1 R(\xi, X)\xi-\varphi R(\xi,\varphi X)\xi=0.
\end{eqnarray}
\end{cor}

Below are the conclusions of Theorem \ref{opkrzywnorm}. concerning the Ricci curvature.

\begin{cor}
The Ricci Tensor of the conformally flat normal almost (para)\-con\-tact metric manifold of dimension $2n+1\geqslant5$ satisfies the relation
\begin{eqnarray}
\label{ric0}
  && \widehat{\mathop{Ric}}\,\varphi X-\varphi\widehat{\mathop{Ric}}\,X
     =\eta(\widehat{\mathop{Ric}}\,\varphi X)\xi-\eta(X)\varphi\widehat{\mathop{Ric}}\,\xi,\\
\label{ric1}
  && \mathop{Ric}(X,Y)+\varepsilon_1\mathop{Ric}(\varphi X,\varphi Y)
     =\eta(X)\mathop{Ric}(Y,\xi)+\eta(Y)\mathop{Ric}(X,\xi)\\
  &&\hskip 5.5cm\null-\eta(X)\eta(Y)\mathop{Ric}(\xi,\xi),\nonumber\\
\label{ric2}
  && \mathop{Ric}(X,X)+\varepsilon_1\mathop{Ric}(\varphi X,\varphi X)=0\quad 
     {\it for}\quad X,Y\in \mathcal{D}.
\end{eqnarray}
\end{cor}

\begin{pf}
The Riemannian curvature tensor of the conformally flat metric manifold has the form 
\begin{eqnarray}
\label{flat2}
  R(X,Y)Z&=&\frac{1}{2n-1}(g(Y,Z)\widehat{\mathop{Ric}}\,X+\mathop{Ric}(Y,Z)X\\
         & &\null-g(X,Z)\widehat{\mathop{Ric}}\,Y-\mathop{Ric}(X,Z)Y) \nonumber\\
         & &\null-\frac{r}{2n(2n-1)}(g(Y,Z)X-g(X,Z)Y). \nonumber
\end{eqnarray}
From (\ref{nr9}), after using (\ref{flat2}), (\ref{eps1}) and (\ref{eps2}), we obtain  
$$
  \widehat{\mathop{Ric}}\,\varphi^2X-\varphi\widehat{\mathop{Ric}}\,\varphi X
  =\eta(\widehat{\mathop{Ric}}\,\varphi^2X)\xi.
$$
Putting in the above equality $\varphi X$ instead of $X$ and using once again (\ref{eps1}), we get (\ref{ric0}). Projecting  (\ref{ric0}) onto $\varphi Y$ and using (\ref{eps1}) and (\ref{eps2}), we obtain the condition  (\ref{ric1}). The equality (\ref{ric2}) is a direct consequence of the formula (\ref{ric1}).
\end{pf}

\end {document}